\documentclass[12pt]{article}

\newtheorem{mytheo}{Theorem}[section]
\newtheorem{myprop}{Proposition}[section]

\newtheorem{mylemma}{Lemma}[section]

\begin{document}

\title {\bf Las Vergnas cube conjecture and reconstruction properties of the cube matroid}
\maketitle

\vskip 1cm 
\begin{center}

{\Large \bf Ilda P. F. Da Silva}
\vskip 1cm 

CELC/University of Lisbon

Faculdade de Ci\^encias - Dep.Matem\'atica

Edif\'icio C6 - Piso 2

P - 1749-016 - Lisbon, PORTUGAL

{\bf E-mail:} isilva@cii.fc.ul.pt

{\bf Fax:} 351-21-7500072

{\bf Phone:} 351-21- 7500330
\end{center}

\vfill\eject
.\vskip 3cm
\begin{center} 
{\Large\bf Abstract}
\end{center}

Las Vergnas Cube Conjecture states that the cube matroid  has exactly one class of orientations. We prove that this conjecture is equivalent to saying that the oriented matroid  $Aff(C^n)$, of the affine dependencies of  the n-cube  $C^n:=\{-1,1\}^n$  over  $I\!R$,   can be reconstructed from the underlying matroid  and one of the following partial lists of signed circuits or cocircuits: 1) the signed circuits of rank 3  or    2) the positive signed cocircuits. 
\vskip 2cm 

{\parindent=0cm
{\bf Keywords:} Cube matroid, oriented matroid, orientation class, reconstruction.\vskip 3mm

{\bf AMS classification:} 52C40, 05B35; 52B12, 52B40, 51M20, 15A35 }

\vskip 3cm

\vfill\eject

\section {Introduction}

There are matroids which are not orientable. There are matroids which have more then one class of orientations. There are matroids which have exactly one class of orientations.

The general problems of determining the orientability of a matroid and the number of reorientation classes of an orientable matroid were considered in the seminal paper of R. Bland and M. Las Vergnas \cite {BLV} where,in particular, it was  shown that regular matroids have exactly one class of orientations. 
 
Asymptotic bounds for the number of orientation classes of realizable uniform oriented matroids can be obtained from \cite {GP}, \cite {A}  and for uniform oriented matroids (not necessarily realizable) are given in \cite {OM}. 

In this paper we consider the following conjecture of M. Las Vergnas concerning the number of orientation classes of  the {\it cube matroid}, i.e. the matroid  of the affine dependencies of  $C^n$, the set of vertices of an n-dimensional cube of  $I\!R^n$:
\vskip 2mm

{\parindent=0cm
{\bf Las Vergnas Cube Conjecture:} \cite {LV}, \cite {OM}, \cite {B}
{\it The cube matroid has a unique class of orientations}.}
\vskip 2mm

The conjecture was proven to be true for  $n\leq 4$  by M. Las Vergnas, J.-P. Roudneff and I. Sala\"un  in \cite {LV}. Later,  J. Bokowski, A. Guedes de Oliveira, U. Thiemann and A. Veloso da Costa \cite {B} verified the conjecture for  $n\leq 7$. 

The main result of this paper is  Theorem 3.1 which states that every class of orientations of the cube matroid has an orientation which coincides on the rank 3 circuits or equivalently on the positive cocircuits with the orientation  $Aff(C^n)$. As a consequence of this theorem we obtain  Theorem 3.2  which reestates Las Vergnas Cube Conjecture in terms of  reconstruction properties of the signatures of circuits and cocircuits of the oriented matroid  $Aff(C^n)$.

The results are presented in the next sections 2 and 3. Section 2  is devoted to properties of the cube (matroid) and Section 3 to properties of the orientations of the cube. Some final remarks are presented in section 4.  

We assume that the reader is familiar with matroid and oriented matroid terminology. Good references are  \cite {OM}, \cite {W}, \cite {M}.

\section {The Cube Matroid}

In this section we develop some terminology and notation to handle the cube matroid. We introduce the notion of  {\it k-subcube}  of the cube matroid (see Definitions 2.1. and Theorem 2.1).  We then present some further properties of the  {\it (n-1)-subcubes}, the facets and skew-facets of the cube matroid and of  the  {\it 2-subcubes or rectangles}  which are the rank 3 circuits of the cube matroid.

We start by recalling that apart from the definition there is no known explicit description of the cube matroid  for every  dimension  $n\in I\!N$. In \cite {IS1} the reader can find an explicit description in terms of hyperplanes for dimension up till  $7$. 
\vskip 5mm

{\parindent=0cm
{\bf Notation. }  

We consider as standard n-cube the set  $C^n:=\{-1,1\}^n$. 

An element of  $C^n$  is called a vertex of the n-cube and is denoted  ${\bf v}=(v_1,\ldots, v_n)$  or simply  $\bf v$.}

Given a vertex  ${\bf v}=(v_1,\ldots,v_n)\in C^n$  and a subset  $I\subseteq  [n]$  we denote by   $_{-I}{\bf v}$  the vertex whose entries are obtained reversing the signs on the entries of  $\bf v$ indexed by $I$  and by  ${\bf v}(I)$  the vector obtained replacing by zeros the entries of  $\bf v$  indexed by  $[n]\setminus I$:
$$_{-I}{\bf v}=(v'_1,\ldots,v'_n)\ \   where\ \ v'_i=-v_i\ \ if \ i\in I \ \ and \ \ v'_i=v_i\ \ if \ \ i\in [n]\setminus I.$$
$${\bf v}(I)=(v'_1,\ldots,v'_n)\ \ where \ \ v'_i=v_i\ \ if \ i\in I\ \  and \ \  v'_i=0\ \ if \ i\in [n]\setminus I.$$
While $_{-I}{\bf v}$  is a new vertex of  $C^n$, ${\bf v}_I$  is as a vector of  $I\!R^n$. The following equality holds:  $_{-I}{\bf v}={\bf v}-2{\bf v}(I)$.  

If  $I,J\subseteq [n] $  are  {\it disjoint } then we write  $I\uplus J\subseteq [n]$. In this case, to simplify notation, we use  $IJ$  instead of  $I\cup J$  in variables depending on subsets of  $[n]$. For instance,  given  ${\bf v}\in C^n$  we write   $_{-IJ}{\bf v}$  instead of   $_{-(I\uplus J)}{\bf v}$  and  if  $J=\{j\}$  we write  $_{-Ij}{\bf v}$.     

{\it The  matroid of affine dependencies of  $C^n$  over  $I\!R$}  will be denoted  $M(C^n)$. We  refer to flats, hyperplanes, cocircuits, circuits,etc... of the matroid  $M(C^n)$  as  flats, hyperplanes, cocircuits, circuits,etc... of the  $n$-cube  $C^n$. 

A {\it hyperplane}  $H$  of the  $n$-cube is a subset  $H\subseteq C^n$  that satisfies the following two conditions:  1) The affine span, $aff(H)$,  of  $H$  is an affine hyperplane of  $I\!R^n$  and  2) $H=aff(H)\cap C^n$.

We identify a hyperplane  $H$  of  $C^n$  with a linear equation defining the affine hyperplane $aff(H)$  of  $I\!R^n$.  When we refer to the hyperplane of  $C^n$  defined by  $H:{\bf x}.{\bf u}=b$, for some fixed  ${\bf u}\in I\!R^n$  and  $b\in I\!R$,  we mean that the hyperplane  $H$  is the set of solutions  ${\bf v}\in C^n$  of this linear equation.

Between the hyperplanes of  the  $n$-cube we distinguish the {\it facets and skew-facets}. Denoting by  $({\bf e_1}, \ldots, {\bf e_n})$  the canonical basis of  $I\!R^n$,  The {\it facets of  $C^n$}  are the  $2n$  hyperplanes  defined by  $H_{\epsilon i}: {\bf x}.{\bf e_i}=\epsilon$, $\forall i\in [n]$  and  $\epsilon\in \{-1,+1\}$. The  {\it skew - facets of  $C^n$}  are the  $n^2+n$ hyperplanes  defined by  $H_{\epsilon ij}: {\bf x}.({\bf e_i}+\epsilon {\bf e_j})=0$, $\forall i<j\in [n]$  and  $\epsilon\in \{-1,+1\}$.   

   We recall that {\it a cocircuit of  $C^n$}  is the complement  $C^n\setminus  H$  of a hyperplane  $H$  of  $C^n$. {\it A circuit of  $C^n$}  is a subset  $C=\{ {\bf v_1,\ldots,v_{r+1}} \}$  of vertices of  $C^n$  which is minimal affine dependent, i.e.  $C$  is affinely dependent  and  $C\setminus  \{\bf v_i\}$  is affinely independent  $\forall  {\bf v_i}\in C$.

{\it The rank, $r(A)$, of a subset  $A\subseteq C^n$}  is related to the dimension of the affine span of  $A$  in the following way:  $r(A)=dim(aff(A))+1$. In particular, the rank of a circuit with  $r+1$  elements is  $r$.

\subsection {Subcubes of  $C^n$}

{\parindent=0cm
{\bf Definitions 2.1.1.} A {\it k-subcube of  $C^n$}  is a subset  $C\subset C^n$  such that the matroid of affine dependencies of  $C$  over  $I\!R$  is isomorphic to the matroid of the  {\it k-cube}  $C^k$.

{\bf 2.1.2.} The  {\it $k$-subcube of  $C^n$  generated by a vertex  $\bf v$  and a  $k$-partition  $I_1\uplus I_2\uplus\ldots \uplus I_k$  of a subset of  $[n]$}, denoted  $C ({\bf v}; I_1,\ldots,I_k)$,  is the subset of  $C^n$  defined by:

$C({\bf v};I_1,\ldots,I_k):=\{ {\bf w}\in C^n:\ {\bf w}=\epsilon_1{\bf v}(I_1)+\ldots +\epsilon_k{\bf v}(I_k)+{\bf v}(J)$  where  $\ J=[n]\setminus (I_1\uplus \ldots \uplus I_k)\ \ and\ \ \epsilon_i\in \{-1,1\} \}$. } 
\vskip 2mm

The next Theorem characterizes the subcubes of  $C^n$. 
\vskip 5mm

{\parindent=0cm
{\mytheo {For a subset  $C\subseteq  C^n$  the following four conditions are equivalent:
\begin{enumerate}
\item $C$  is a  $k$-subcube of  $C^n$.
\item $|C|=2^k$  and  $dim(aff(C))=k$.
\item $C$  is a flat of  $C^n$  with rank  $k+1$  and maximum number of elements.
\item $C$  is the k-subcube  $C({\bf v};I_1,\ldots,I_k)$   generated by a vertex  ${\bf v}\in C^n$  and some  $k$-partition  $I_1\uplus \ldots \uplus I_k$  of a subset of  $[n]$.
\end{enumerate}
}}

The proof of this theorem consists in showing the following implications:  $1)\Longrightarrow 2) \Longleftrightarrow 3) \Longrightarrow 4) \Longrightarrow 1)$. The implications which are not obvious are  $2) \Longleftrightarrow 3)$  and  $2),3)\Longrightarrow 4)$. They are proved by double induction: first on  $k$  then on  $n$. }

We would like to mention that a version of the next Lemma 2.1. appears in \cite{B}.
\vskip 5mm

{\parindent=0cm
{\mylemma{ For a subset  $H\subseteq C^n$  the following four conditions are equivalent:
\begin{enumerate}
\item $H$  is a  $(n-1)$-subcube of  $C^n$.
\item $|H|=2^{n-1}$  and  $dim(aff(H))=n-1$.
\item $H$  is a hyperplane of  $C^n$  with maximum number of elements.
\item $H$  is a facet or a skew facet of  $C^n$.
\end{enumerate}
}}
\vskip 2mm

{\bf Proof.}} The implications  $1)\Longrightarrow 2)$  and  $4)\Longrightarrow 1)$  are obvious. We prove by induction on  $n$  that  $2)\Longleftrightarrow 3)$  and  $2),3)\Longrightarrow 4)$.

The equivalence  $2)\Longleftrightarrow 3)$  is a direct consequence of the following claim:

{\parindent=0cm
{\bf Claim 1:} {\it If  $G\subseteq C^n$  is such that  $|G|\geq 2^{n-1}+1$  then  $dim(aff(G))=n$.}}

The proof of this claim is by induction on  $n$. The claim is clearly true for  $n=1,2,3$. Assume the claim is true for  $n$  and consider  $G\subseteq C^{n+1}$  such that   $|G|\geq 2^n+1$. Consider the facets  $H_1,\ H_{-1}$  of  $C^n$  and define  $G_1:=G\cap H_1$  and  $G_{-1}:=G\cap H_{-1}$. Since  $G=G_1\uplus G_{-1}$  and  $|G|\geq 2^n+1$, $G_1$  and  $G_{-1}$  are both nonempty and one of these sets, say  $G_1$, contains at least  $2^{n-1}+1$  elements. Since  $H_1$  is a  $n$-cube, the induction assumption implies that  $dim (aff(G_1))=n$  and consequently that  $dim(aff(G))=dim (aff (G_1\uplus G_{-1}))=n+1$.

The implication $ 2),3) \Longrightarrow 4)$  is a direct consequence of the next claim:

{\parindent=0cm
{\bf Claim 2:} {\it Let  $H$  be a hyperplane of  $ C^n$  with  $2^{n-1}$  elements. Then, either  $H$  is a facet:  $H_{\epsilon i}:{\bf x. e_i}=\epsilon $  or  $H$  is a skew facet:  $H_{\epsilon ij}:{\bf x.}({\bf e_i}+\epsilon {\bf e_j})=0$, for some $\epsilon\in \{-1,1\}$} }

The proof of this claim is also by induction on  $n$. The claim is clearly true for  $n=1,2,3$. Assume the claim is true for  $n$  and consider a hyperplane  $H$  of  $C^{n+1}$  such that   $|H|=2^n$. Two cases are possible: 

{\parindent=0cm
Case 1)  $|H\cap H_{\epsilon i}|\geq  2^{n-1}+1$, for some  $i\in [n+1],\ \epsilon\in \{-1,1\}$. In this case Claim 1 implies that  $H=H_{\epsilon i}$ and claim 2 follows.

Case 2)} $|H\cap H_{\epsilon i}|=  2^{n-1}$,   $\forall i\in [n+1]\ \epsilon \in\{-1,1\} $.   In this case we consider the facets  $H_1,\ H_{-1}$  of  $C^n$  and let  $G_1:=H\cap H_1$  and  $G_{-1}:=H\cap H_{-1}$. By the induction assumption  either  

{\parindent=0cm
(A)  $G_1: ({\bf x. e_1}=1 \ \ and\ \ {\bf x. e_i}=\epsilon )\ \   for\ some\  \ i\in \{2,\ldots,n+1\},$ $\ \epsilon\in \{-1,1\}$

or
  
(B)  $G_1: ({\bf x. e_1}=1 \  and\ \ {\bf x.}({\bf e_i}+\epsilon {\bf e_j})=0) \ for\ some\  \ i,j\in \{2,\ldots, n+1\},$ $\ i\not=j, \ \epsilon\in \{-1,1\}$.}

In case (A), since  $H=G_1\uplus G_{-1}$  is not a facet of  $C^n$  we must have  $G_{-1}: ({\bf x . e_1}=-1 \ \ and\ \ {\bf x . e_i}=-\epsilon )$  implying that, in this case  $H$  is the skew-facet defined by  $H_{\epsilon 1j}:{\bf x.} ({\bf e_1}+\epsilon {\bf e_j})=0$.

In case (B),  $H$  must be the skew-facet  $H_{\epsilon ij}:{\bf x .}({\bf e_i}+\epsilon {\bf e_j})=0$. In fact, if  $H\not= H_{\epsilon ij}$  since  $|G_1\cap H_{\epsilon ij}|=2^{n-1}$  and  $C^n=H_{\epsilon ij}\uplus H_{-\epsilon ij}$  we would have  
$|G_{-1}\cap H_{-\epsilon ij}|=2^{n-1}$  but this implies that  $dim (aff(H))=n+1$, contradicting the assumption that  $H$  is a hyperplane of  $C^{n+1}$.   
\vskip 5mm 

{\parindent=0cm
{\bf Proof of Theorem 2.1} } Using lemma 2.1 we prove by double induction: first on  $k$  then on  $n$  the non-obvious implications:  $2) \Longleftrightarrow 3)$  and  $2),3)\Longrightarrow 4)$. 

In order to prove the equivalence  $2) \Longleftrightarrow 3)$  we prove the following:

{\parindent=0cm
{\bf Claim 1:} {\it If  $C\subseteq C^n$  is such that  $|C|\geq 2^{k}+1$  then  $dim(aff(C))\geq k+1$  (or equivalently  $r(C)\geq k+2$).} }

This claim is clearly true for  $k=0$  and all $ n\in I\!N$.  

Assume that the claim has been proved for  $k$  and  all  $n\in I\!N,\ n\geq k+1$.

Consider  $C\subseteq C^n$  such that  $|C|\geq 2^{k+1}+1$. Then  $n\geq k+2$. If  $n=k+2$  then Lemma 2.1. implies that  $dim(aff (C))=k+2$  and the claim is verified.

Assume that the claim is true for  $k+1$  and all  $n$  such that  $k+1\leq n<m$  and consider  $C\subseteq C^m$  such that  $|C|\geq 2^{k+1}+1$. If  $C$  is contained in some facet  $H_{\epsilon i}$  of  $C^m$  then, since  $H_{\epsilon i}$  is a  (m-1)-cube, the induction assumption guarantees that the claim is verified. 

We now consider the case  where  $C$  is contained in none of the facets of  $C^m$. 

Define  $C_1:= C\cap H_1$  and  $C_{-1}=C\cap H_{-1}$. Since   $|C|\geq 2^{k+1}+1$  both these sets are nonempty and one of them, say  $C_1$  must  contain more then  $2^{k}+1$  elements. The induction assumption implies then that  $dim(aff(C_1))\geq k+1$  and consequently  $dim(aff ( C_1\cup C_{-1}))=dim (aff(C))\geq k+2$.

The implication  $2),3)\Longrightarrow 4)$  is a direct consequence of the next claim:

{\parindent=0cm
{\bf Claim 2:} {\it Let  $C$  be a flat of  rank  $k+1$ (i.e. $dim(aff(C)=k$)  of  $C^n$  with  $2^{k}$  elements. Then, $C$  is the k-subcube  $C({\bf v}; I_1,\ldots,I_k)$   generated by a vertex  ${\bf v}\in C^n$  and a  $k$-partition  $I_1\uplus \ldots \uplus I_k$  of a subset of  $[n]$.}}

Claim 2 is trivially true for  $k=1$, $\forall n\in I\!N$ since in this case we must have  $C=\{ {\bf v}, _{-I}{\bf v}\}=C({\bf v};I)$  for some vertex  ${\bf v}\in C^n$  and some subset   $I\subseteq [n]$.

Assume that Claim 2 is true for  $k$  and all  $n\in I\!N$  such that  $k\leq n< m$. 

Consider  $C\subseteq C^m$  such that  $|C|=2^{k+1}$. If  $m=k+1$  the result is obvious. If  $m=k+2$ the claim is true by  Lemma 2.1. Assume now that the claim is true for  all $n$  such that  $k+2\leq n <m$  and consider  $C\subseteq C^m$  such that  $|C|=2^{k+1}$  and  $dim(aff(C))=k+1$.

If  $C$  is contained in a facet  $H_{\epsilon i}$  of  $C^m$  then, by the induction assumption there is a partition  $I_1\uplus \ldots \uplus I_{k+1}\uplus J$  of a subset of $[m]\setminus \{i\}$  (eventually  $J=\emptyset$) such that  being  $\bf v$  a vertex of  $C$,  $C\cap H_{\epsilon i}$  is the  $(k+1)-cube$  of  $H_{\epsilon i}$  generated by the partition  $I_1\uplus \ldots \uplus I_{k+1}\subseteq [m]\setminus \{i\}$  which is the cube  $C=C({\bf v};I_1,\ldots,I_{k+1}):=\{ \epsilon_1{\bf v}({I_1})+\ldots +\epsilon_k{\bf v}({I_k})+{\bf v}({Ji})\ where \  \epsilon_i\in \{-1,1\} \}$  of  $C^m$. The claim is verified in this case.

If  $C$  is not contained in a facet of  $C^m$  then, by Claim 1  each facet of  $C^m$  contains  exactly  $2^{k}$  points of  $C$. Consider  $C_1:=C\cap H_1$  and  $C_{-1}:=C\cap H_{-1}$.  By the induction assumption there is a  $k$-partition  $I_1\uplus \ldots \uplus I_k$  of a subset of  $[m]\setminus \{1\}$  such that  $C_1$  is a  $k-subcube$  of  $H_1$  i.e.   $C_1:={\bf v}({J})+\epsilon _1{\bf v}(I_1)+\ldots +\epsilon_k {\bf v}({I_k})$, for some vertex  $\bf v\in C_1$, some  k-partition  $I_1\uplus\ldots,\uplus I_k\subseteq [m]\setminus \{1\}$ and  $J=[m]\setminus I_1\uplus\ldots,\uplus I_k$. This implies that  $aff(C_1):={\bf v}({J})+ lin ({\bf v}({I_1}),\ldots , {\bf v}({I_k}))$. Note that  $J\not=\emptyset$  because  $1\in J$.

On the other hand,  $aff(C_{-1})$  is an affine subspace of  $I\!R^n$  paralell to  $aff(C_1)$  and with the same dimension therefore, for any  ${\bf w}\in C_{-1}$  we have:  $aff(C_{-1}):={\bf w}({J})+ lin ({\bf v}({I_1}),\ldots , {\bf v}({I_k}))$. Since  $C=C_1\uplus C_{-1}$  is not contained in  $H_{\epsilon i},\ \forall i\in [m]$  we must have  ${\bf v}({J})=-{\bf w}({J})$, implying that  $C$  is the  (k+1)-cube  $C=C({\bf v};I_1,\ldots,I_{k},J)$.
\vskip 5mm

{\parindent=0cm
{\bf Remark 2.1.} It is clear from Lemma 2.1. that the  (n-1)-cubes of  $C^n$  are the facets and skew-facets, the hyperplanes of the matroid  $M(C^n)$  with largest number of elements. 

Theorem 2.1. shows, in particular, that the  2-subcubes of  $C^n$  are the subsets of the form  $C=C({\bf v}; I,J)=\{ {\bf v}, _{-I}{\bf v}, _{-IJ}{\bf v}, _{-J}{\bf v}\}$  with  $I\uplus J\subseteq [n]$  i.e. the four vertices of a rectangle of  $I\!R^n$  and therefore we call them {\it rectangles of  $C^n$}.  The next proposition shows that the rectangles of  $C^n$  are precisely the rank 3 circuits of the cube matroid.  }
\vskip 5mm

{\parindent=0cm
{\myprop {Consider three distinct vertices  ${\bf v,v_1,v_2}\in C^n$  and the affine plane   $P:=aff({\bf v,v_1,v_2})$. Let  $I,J$  be the subsets of  $\left[ n\right]$  such that  ${\bf v_1}=_{-I}{\bf v}$  and  ${\bf v_2}=_{-J}{\bf v}$. Then one (and only one) of the following three situations must occur:
\begin{enumerate}
\item $I\cap J=\emptyset $  and  $P\cap C^n=\{ {\bf v,v_1,v_2,v_3} \}$  where  ${\bf v_3}=_{-IJ}{\bf v}$.
\item $I\subset J $ (or  $J\subset I$)  and  $P\cap C^n=\{ {\bf v,v_1,v_2,v_3} \}$  where  ${\bf v_3}=_{-(J\setminus I)}{\bf v}$  (resp. ${\bf v_3}=_{-(I\setminus J)}{\bf v}$).
\item $I\cap J, I\setminus J$  and  $J\setminus I$  are nonempty and  $P\cap C^n=\{ {\bf v,v_1,v_2,} \}$.
\end{enumerate}}}

{\parindent=0cm
{\bf Proof.}} The plane  $P$  is the set of all the affine combinations:
$${\bf p}_{a,b}=(1-a-b){\bf v}+a_{-I}{\bf v}+b_{-J}{\bf v},\ \ a,b\in I\!R$$
of the points  ${\bf v,v_1}=_{-I}{\bf v},v_2=_{-J}{\bf v} $.

If  ${\bf v}=(v_1,\ldots,v_n)\in C^n$. Then the i-th coordinate,  $({\bf p}_{a,b})_i$, of  ${\bf p}_{a,b}$  is given by:
$$({\bf p}_{a,b})_i=\left\{ \begin{array}{cc}(1-2a-2b)v_i&\ i\in I\cap J \\  (1-2a)v_i& \ i\in I\setminus J\\ (1-2b)v_i&\ i\in J\setminus I\\  v_i&\ i\in \left[  n\right] \setminus (I\cup J)\end{array}\right.$$

1) If  $I\cap J=\emptyset$  then  ${\bf p}_{a,b}\in C^n$  iff  $1-2a$,$1-2b\in \{-1,1\}\Longleftrightarrow \ a,b\in \{0,1\}$. In this case  $P\cap C^n=\{{\bf p}_{0,0}={\bf v},{\bf p}_{1,0}={\bf v_1},{\bf p}_{0,1}={\bf v_2}, {\bf p}_{1,1}= _{-IJ} {\bf v}\}$.

2) If  $I\subset J$  then  ${\bf p}_{a,b}\in C^n$  iff  $1-2a-2b$,$1-2b\in \{-1,1\}$. Therefore either   $b=0$  and  $a\in \{0,1\}$  or  $b=1$  and  ${\bf a}\in \{-1,0\}$. In this case  $P\cap C^n=\{{\bf p}_{0,0}={\bf v},{\bf p}_{1,0}={\bf v_1},{\bf p}_{1,0}={\bf v_2}, {\bf p}_{-1,1}= _{J\setminus I} {\bf v}\}$.

3) If  $I\cap J, I\setminus J$  and  $J\setminus I$  are non empty then  ${\bf p}_{a,b}\in C^n$  iff  $1-2a-2b$, $1-2a$,$1-2b\in \{-1,1\}$. There are 3 pairs  $(a,b)$  satisfying these conditions the pairs  $(0,0)$, $(1,0)$  and  $(0,1)$  implying that, in this case,  $C^n\cap P= \{{\bf p}_{0,0},{\bf p}_{1,0},{\bf p}_{0,1}\}$ $={\{\bf v,v_1,v_2}\}$.

\vskip 5mm

The next proposition translates in terms of rectangles  the elimination property for modular pair of circuits. The proof is left to the reader.
\vskip 5mm 

{\parindent=0cm
\myprop {Let  ${\bf v}=(v_1,\ldots,v_n)\in C^n$  and  $I,J,K$  be three disjoint nonempty subsets of  $[n]$. 

\begin{enumerate}
\item  Consider the rectangles  $C({\bf v};I,J)$  and  $C({\bf v};I,K)$  then the unique circuit of the cube matroid contained  in  $(C({\bf v};I,J) \cup C({\bf v};I,K))\setminus  \{\bf v\}$  is the circuit   $C(_{-J}{\bf v};I,JK)$. 

\item Consider the rectangles  $C({\bf v};IJ,K)$  and  $C({\bf v};I,JK)$  then the unique circuit of the cube matroid contained  in  $(C({\bf v};IJ,K) \cup C({\bf v};I,JK))\setminus  \{\bf v\}$  is the circuit   $C(_{-I}{\bf v};J,IK)$.  
\end{enumerate}}}
}

\vskip 5mm

\section {Orientations of the Cube Matroid}
 
\subsection { The oriented matroid  $Aff(C^n)$}

The oriented matroid of affine dependencies of $C^n$  over  $I\!R$, denoted  $Aff(C^n)$, is the  orientation of the cube matroid  $M(C^n)$  whose {\it signature of cocircuits  $\mathcal D$} is defined in the following way:

Consider a  cocircuit  $Y$  of  $C^n$. The complement  $H=C^n\setminus Y$  of  $Y$  is a hyperplane $H: {\bf x . h}=b$  of  $C^n$. Consider the partition of  $Y$  into the subsets   $Y^+:=\{ {\bf v}\in C^n: \ {\bf v .}{\bf h}>b\}$  and  $Y^-:=\{ {\bf v}\in C^n: \ {\bf v .}{\bf h}<b\}$. The signature of the cocircuit  $Y$  in the orientation  $Aff(C^n)$  is the pair of opposite signed sets  $Y=(Y^+,Y^-)$  and  $-Y=( Y^-, Y^+)$. 
 
Note that the positive cocircuits of  $Aff(C^n)$  are the cocircuits  $Y_{\epsilon i}:= (H_{\epsilon i}, \emptyset)$  complementary of the hyperplanes   $H_{\epsilon i}:{\bf x .}{\bf e_i}=\epsilon$, {\it the facets of  $C^n$}.

We denote by  ${\mathcal F}$  the subfamily of  $\mathcal D$  which contains the positive cocircuits and its opposites:
$${\mathcal F}:=\{ \pm Y_{\epsilon i}\in \mathcal D: \ Y_{\epsilon i}=(H_{-\epsilon i},\emptyset),\ i\in [n],\ \epsilon \in \{-1,1\} \}$$

The {\it signature of circuits  $\mathcal C$}  of the orientation  $Aff(C^n)$  is defined in the following way:

Given a circuit  $X$  of  $C^n$  there is a unique partition of  $X$  into two disjoint subsets  $X=X^+\uplus X^-$  with the property that  $conv(X^+)\cap conv (X^-)\not=\emptyset$. The signature of the circuit  $X$  is the pair of opposite signed sets  $X=(X^+,X^-)$, $-X=(X^-,X^+)$.

The rank three circuits of  $Aff(C^n)$, the {\it signed rectangles}, are the signed subsets of the form  $\pm R ({\bf v}; {I,J})$  with  $I\uplus J\subseteq [n]$ defined by:
$R=R ({\bf v}; {I,J})= (\{{\bf v}, _{-IJ}{\bf v}\},  \{_{-I}{\bf v}, _{-J}{\bf v}\})=  {\bf v}^+\  _{-I}{\bf v}^-  \ _{-IJ}{\bf v}^+\  _{-J}{\bf v}^-$ 

We denote by  
$\mathcal R$  the subfamily of  $\mathcal C$  which contains the signed rectangles of  $C^n$:
$${\mathcal R}:=\{ \pm R({\bf v};I,J)= {\bf v}^-\  _{-I}{\bf v}^+  \ _{-IJ}{\bf v}^-\  _{-J}{\bf v}^+,\  I\uplus J\subseteq [n], {\bf v}\in C^n\}.$$
\vskip 2mm

We recall that the families of signed circuits and cocircuits of an oriented matroid are orthogonal. In what follows we will make extensive use of this property  which we briefly recall:

{\it Two signed subsets  $X=(X^+,X^-), Y=(Y^+,Y^-)$  of a set  $E$  are orthogonal}, written  $X\perp Y$  iff the following condition is satisfied:
$$(X^+\cap Y^+)\cup (X^-\cap Y^-)\not= \emptyset\ \ iff\ \ (X^+\cap Y^-)\cup (X^-\cap Y^+)\not= \emptyset \leqno(O)$$
{\it Two families  $\mathcal X$  and  $\mathcal Y$  of signed subsets of  $E$  are orthogonal} if  $\forall X\in {\mathcal X}, Y\in {\mathcal Y}\ \ \ X\perp Y$. For more details see \cite {BLV}, \cite {OM}.

\vskip 2mm

{\parindent=0cm
{\bf Remark 3.1.}  In what follows  $\mathcal C$  and  $\mathcal D$  allways represent the signatures, respectively of the circuits and cocircuits of the oriented matroid  $Aff(C^n)$.

${\mathcal F}$  denotes the subfamily of  $\mathcal D$  containing  the positive and negative cocircuits of  $Aff(C^n)$. ${\mathcal R}$  denotes the subfamily of  $\mathcal C$  corresponding to signed rectangles of  $Aff(C^n)$.}

\vskip 5mm

\subsection {Properties of the orientations of the n-cube}

{\parindent=0cm 

\myprop{
For an orientation  ${\mathcal M}$  of the cube matroid  $M(C^n)$   with signatures of cocircuits  and circuits, respectively,  ${\mathcal D}'$  and  ${\mathcal C}'$, the following conditions are equivalent:

\begin{enumerate}
\item  ${\mathcal F}\subseteq {\mathcal D}'$.
\item  ${\mathcal R}\subseteq {\mathcal C}'$.
\end{enumerate}}
\vskip 2mm

{\bf Proof.} } Immediate consequence of the orthogonality between the signatures of circuits and cocircuits of an oriented matroid.
\vskip 2mm 

{\parindent=0cm
\mytheo { For every orientation   ${\mathcal  M}$   of the cube matroid  $M(C^n)$ there is a subset  $A\subseteq C^n$  such that the reorientation  $_{-A}{\mathcal M}$,  obtained from  $\mathcal M$  reversing signs on the subset  $A$,  satisfies one of the following (equivalent) conditions: 
  
\begin{enumerate}
\item  ${\mathcal F}\subseteq {\mathcal D}'$.
\item  ${\mathcal R}\subseteq {\mathcal C}'$.
\end{enumerate}}

where  ${\mathcal C}'$  and  ${\mathcal D}'$  are the families of signed circuits and signed cocircuits of the reorientation $_{-A}{\mathcal M}$.
\vskip 2mm
}
\vskip 2mm

The proof of this theorem is done in several steps presented in the next three Lemmas. 
\vskip 5mm

{\parindent=0cm
\mylemma{ Let  ${\mathcal M}$  be an orientation of the cube matroid  $M(C^n)$. Then, there is  $B\subseteq C^n$  such that  $_{-B}{\mathcal M}$  is acyclic and contains the positive cocircuits  $Y_{-n}=(H_{n},\emptyset)$  and  $Y_n= (H_{-n}, \emptyset)$.}
\vskip 2mm

{\bf Proof.}} Consider an orientation  $\mathcal M$  of the n-cube matroid. Let  $X_n, X_{-n}$  denote the signed cocircuits of  $\mathcal M$  complementary of the hyperplanes  $H_n$  and  $H_{-n}$ , respectively. Then   $X_n=(H^{+}_{-n}, H^{-}_{-n})$  for some partition  $H^{+}_{-n}\uplus H^{-}_{-n}$  of  $H_{-n}$  and  $X_{-n}=(H^+_n, H^-_n)$  for some partition  $H^{+}_n\uplus H^{-}_n$  of  $H_n$.  Define  $B:= H^{-}_{-n}\cup H^{-}_n$. The reorientation  $_{-B}{\mathcal M}$  contains the positive cocircuits  $Y_{n}=_{-B}X_{n}=(H_{-n},\emptyset )$  and   $Y_{-n}=_{-B}X_{-n}=(H_n,\emptyset )$. Since  $C^n=H_n\cup H_{-n}$  we conclude that   $_{-B}{\mathcal M}$  is acyclic.
\vskip 5mm

{\parindent=0cm
{\bf Remark  3.2.} If  $\mathcal M$  is an acyclic orientation of the cube matroid  $M(C^n)$  containing  the positive cocircuits  $Y_n=(H_{-n},\emptyset)$  and  $Y_{-n}=(H_{n},\emptyset)$  then by orthogonality with these cocircuits  the signed  rectangles of  $\mathcal M$  whose support is a rectangle of the form  $C({\bf v};I,Jn)$  with  $I\uplus Jn\subseteq \left[ n \right] $  are 
either  $\pm R({\bf v};(I,Jn))$  or   $\pm R'({\bf v};(I,Jn))$  where  
$$R({\bf v};I,J)={\bf v}^+\ _{-I}{\bf v}^-\ _ {-IJn}{\bf v}^+ \ _{-Jn}{\bf v}^-\ \ \ and\ \ \ R'({\bf v};I,J)={\bf v}^+\ _{-I}{\bf v}^-\ _{-IJn}{\bf v}^-\ _{-Jn}{\bf v}^+.$$  

\vskip 5mm

\mylemma{ Let  ${\mathcal M}={\mathcal M}(C^n)$  be an acyclic orientation of  $C^n$  containing the positive cocircuits  $Y_{-n}=(H_{n},\emptyset)$  and  $Y_n= (H_{-n}, \emptyset)$, then  $\mathcal M$  satisfies one (and only one) of the following conditions:
\begin{enumerate}
\item For every  ${\bf v}\in C^n$  and every  $I\subseteq  \left[ n-1 \right] $, $I\not=\emptyset$  $$R({\bf v};I,n)={\bf v}^+\	 _{-I}{\bf v}^-\ _{-In}{\bf v}^+\ _{-n}{\bf v}^-$$  is a signed circuit of  $\mathcal M$.
\item For every  ${\bf v}\in C^n$  and every  $I\subseteq  \left[ n-1 \right] $, $I\not=\emptyset$ $$R'({\bf v};I,n)={\bf v}^+\ _{-I}{\bf v}^-\ _{-In}{\bf v}^-\ _{-n}{\bf v}^+$$  is a signed circuit of  $\mathcal M$.
\end{enumerate}
}
\vskip 2mm

{\bf Proof.}} Since a rectangle  $C({\bf v};I,n)$  contains  ${\bf v}$  and  $_{-n}\bf v$ we will assume,without loss of generality, that  ${\bf v}\in H_n$.
  
First we prove that {\it for a fixed vertex  ${\bf v}\in H_n$  either $\forall I\subseteq  \left[ n-1 \right] $, $I\not=\emptyset$, $R({\bf v};I,n)$  is a signed circuit of  $\mathcal M$  or $\forall I\subseteq  \left[ n-1 \right] $, $I\not= \emptyset$, $R'({\bf v};I,n)$  is a signed circuit of  $\mathcal M$}. 

Assume, on the contrary, that there are subsets  $I,J\subseteq  \left[ n-1 \right] $, $I, J\not=\emptyset$  such that  $R=R({\bf v};I,n)$  and  $R'=R'({\bf v};J,n)$  are signed circuits of  $\mathcal M$. We consider separately the cases  $I\cap J=\emptyset$  and  $I\cap J\not=\emptyset$  

{\parindent=0cm
{\it Case 1)} If  $I\cap J=\emptyset$  then by Proposition 2.2  we know that there is unique circuit contained in   $(C({\bf v};I,n))\cup C({\bf v};J,n))\setminus \{\bf v\}$  which  is the circuit  $C(_{-I}{\bf v};IJ,n)$. By the elimintaion property for signed circuits of an oriented matroid, the signature of this circuit in  $\mathcal M$,  obtained eliminating  $\bf v$  between  the signed rectangles  $R={\bf v}^+\ _{-I}{\bf v}^-\ _{-In}{\bf v}^+ \ _{-n}{\bf v}^-$  and   $-R'={\bf v}^-\	 _{-J}{\bf v}^+\ _{-Jn}{\bf v}^+ \ _{-n}{\bf v}^-$   must be  $\pm X$  with   $X={_{I}{\bf v}^-\	 _{-In}{\bf v}^+}\ _{-J}{\bf v}^+\ _{-Jn}{\bf v}^+ $. This signed set is not orthogonal to the positive cocircuit  $Y_n$, a contradiction. 

{\it Case 2)} If $I\cap J\not= \emptyset$  then consider  $K:=I\cap J$,  $I_1=I\setminus K$  and  $J_1=J\setminus K$  and the three rectangles:  $C(_{-K}{\bf v};I_1,n)$, $C(_{-K}{\bf v};J_1,n)$   and  $C( _{-K}{\bf v};K,n)$. By the previous case  the signature of circuits of  $\mathcal M$  satisfies one (and only one) of the following two conditions:

A)  The three signed rectangles  $R( _{-K}{\bf v};I_1,n)$, $R( _{-K}{\bf v};J_1,n)$  and  $R( _{-K}{\bf v};K,n)$,  are signed circuits of  $\mathcal M$.

B)  The three signed rectangles  $R'(_{-K}{\bf v};I_1,n)$, $R'( _{-K}{\bf v};J_1,n)$  and  $R'(_{-K}{\bf v};K,n)$,  are signed circuits of  $\mathcal M$.}

If  $\mathcal M$  satisfies condition  $A)$  then eliminating  $_{-K}\bf v$  between  $R({_{-K}{\bf v}};J_1,n)$  and  $-R({_{-K}{\bf v}};K,n)$  we conclude  that  $R({\bf v};J,n)$   must be a signed circuit of $\mathcal M$  contradicting the assumption that  $R'({\bf v};J,n)$  is a signed circuit.   
  
If  $\mathcal M$  satisfies condition  $B)$, eliminating  $_{-K}\bf v$  between  $R'({_{-K}{\bf v}};I_1,n)$  and  $-R'({_{-K}{\bf v}};K,n)$  we conclude  $R'({\bf v};I,n)$  must be a signed circuit of  $\mathcal M$, contradicting the assumption that  $R({\bf v};I,n)$  is a signed circuit.   

To conclude the proof of the lemma  we need to prove that {\it if  ${\bf v}\in H_n$  is such that  $\forall I\subseteq  \left[ n-1 \right] $, $I\not=\emptyset$, $R({\bf v};I,n)$  (resp.  $R'({\bf v};I,n)$  )  is a signed circuit of  $\mathcal M$  then for every  ${\bf w}\in H_n$  also   $R({\bf w}; I,n)$ (resp. $R'({\bf w};I,n)$  )   is a signed circuit of  $\mathcal M$. } 

Assume that  ${\bf v}\in H_n$  is such that  $\forall I\subseteq  \left[ n-1 \right] $, $I\not=\emptyset$, $R({\bf v};I,n)$ ( resp. $R'({\bf v};I,n)$ )  is a signed circuit of  $\mathcal M$. Consider  ${\bf w}\in H_n$. Then  ${\bf w}=_{-I}{\bf v}$  for some  $I\subseteq  \left[ n-1 \right] $, $I\not=\emptyset$  and  $R({\bf v};I,n)$  is signed circuit of  $\mathcal M$. Since  $R({\bf v};I,n)=R({\bf w};I,n)$  we conclude  that  $\forall I\subseteq  \left[ n-1 \right] $, $I\not=\emptyset$, $R({\bf w};I,n)$  (resp. $R'({\bf w};I,n)$)    is a signed circuit of  $\mathcal M$.    
\vskip 5mm

{\parindent=0cm
\mylemma{ Let  ${\mathcal M}={\mathcal M}(C^n)$  be an acyclic orientation of  $C^n$   containing the positive cocircuits  $Y_{-n}=(H_{n},\emptyset)$  and  $Y_n= (H_{-n}, \emptyset)$. Then  $\mathcal M$  satisfies one (and only one) of the following properties:
\begin{enumerate}
\item For every  ${\bf v}\in C^n$  and every  2-partition  $I\uplus J\subseteq  \left[ n-1 \right]  $  of a subset of  $[n-1]$  the signed set:
$$R({\bf v};I,Jn)={\bf v}^+\	 _{-I}{\bf v}^-\ _{-IJn}{\bf v}^+\ _{-Jn}{\bf v}^-$$  is a signed circuit of  $\mathcal M$.
\item For every  ${\bf v}\in C^n$  and every  2-partition  $I\uplus J\subseteq  \left[ n-1 \right]  $  of a subset of  $[n-1]$  the signed set:
$$R'({\bf v};I,Jn)={\bf v}^+\ _{-I}{\bf v}^-\ _{-IJn}{\bf v}^-\ _{-Jn}{\bf v}^+$$  is a signed circuit of  $\mathcal M$.
\end{enumerate}
Moreover, if  $\mathcal M$ satisfies condition  2)  then the orientation, $_{-H_n}{\mathcal M}$  obtained from  $\mathcal M$  reversing signs on  $H_n$  satisfies condition 1).

}
\vskip 2mm

{\bf Proof.}} By Lemma 3.2 we know that  $\mathcal M$  satisfies one (and only one) of the following conditions: 

{\parindent=0cm
A) For every  ${\bf v}\in C^n$  and every  $I\subseteq  \left[ n-1 \right] $, $I\not=\emptyset$
  
$R ({\bf v};I,n)={{\bf v}^+\	 _{-I}{\bf v}^-\ _{-In}{\bf v}^+\ _{-n}{\bf v}^-}$  is a signed circuit of  $\mathcal M$.

B) For every  ${\bf v}\in C^n$  and every  $I\subseteq  \left[ n-1 \right] $, $I\not=\emptyset$  

$R'({\bf v};I,n)={{\bf v}^+\	 _{-I}{\bf v}^-\ _{-In}{\bf v}^-\ _{-n}{\bf v}^+}$  is a signed circuit of  $\mathcal M$.}

We prove that if  $\mathcal M$  satisfies condition A)  then  $\mathcal M$  satisfies condition 1) of the Lemma.

Assume that  $\mathcal M$  satisfies condition A)  and consider  $I\uplus J\subseteq  \left[ n-1 \right] $. Then the signed set  $R({\bf v};IJ,n)$  is a signed circuit of  $\mathcal M$  and by Remark 3.2. either  $R({\bf v}; I,Jn)$  or   $R'({\bf v};I,Jn)$  is a signed circuit of  $\mathcal M$.

Now, if  $R({\bf v};IJ,n)$  and   $R'({\bf v};I,Jn)$  are signed circuits of  $\mathcal M$  Proposition 2.2 implies that  the unique signed circuit of  $\mathcal M$  obtained elimating  $_{-IJn}{\bf v}$  between  $R({\bf v};IJ,n)$  and  $R'({\bf v};I,Jn)$  must be:  

$X=  _{-n}{\bf v}^-\ _{-Jn}{\bf v}^+\ _{-IJ}{\bf v}^-\ _{-I}{\bf v}^-$  which is not orthogonal to the positive cocircuit  $Y_{-n}$. Therefore if  $\mathcal M$  satisfies condition  A)  then  $\mathcal M$  satisfies condition 1).

If  $\mathcal M$  satisfies condition   B)  then  it is clear that  $_{-H_n}{\mathcal M}$  satisfies condition  $A)$  and therefore condition  1)  implying, by the previous case that  $\mathcal M$  satisfies condition 2) of the lemma.
\vskip 5mm

{\parindent=0cm

{\bf Proof of Theorem 3.1.}} Let  $\mathcal M$  be an orientation of the cube matroid  $M(C^n)$.  Consider a subset  $A\subseteq [n]$  such that the reorientation  $_{-A}{\mathcal M}$  satisfies the conditions of Lemma 2.3  i.e.   $_{-A}{\mathcal M}$ contains as signed circuits all  the signed rectangles  $R({\bf v};I,Jn)=$ ${\bf v}^+\ _{-I}{\bf v}^- \ _{-IJn}{\bf v}^+$
 $\ _{-Jn}{\bf v}^-$  with  ${\bf v}\in C^n$  and  $I\uplus J\uplus \{n\}\subseteq [n]$ and as signed cocircuits the positive cocircuits   $Y_n=(H_{-n},\emptyset)$  and  $Y_{-n}=(H_n,\emptyset)$.

We claim that for all  $i\in[n]$ the positive signed sets  $Y_{-i}=(H_i,\emptyset)$  and $Y_{i}=(H_{-i},\emptyset)$  are signed cocircuits of  $_{-A}{\mathcal M}$.

Consider a vertex  ${\bf v}\in H_i$  and let  $X_{-i}=(X_{-i}^+,X_{-i}^-)$  denote the signed cocircuit of  $_{-A}{\cal M}$  with support  $H_{i}$  ($H_i= X_{-i}^+\cup X_{-i}^-$)  such that  ${\bf v}\in X_{-i}^+$. Let  $ _{-I}{\bf v}$  be another vertex of  $H_i$ ,with  $I\subseteq [n]\setminus \{i\}$.  

If  $n\notin I$  then  the signed circuit  $R({\bf v};I,n)={\bf v}^+\ _{-I}{\bf v}^- \ _{-In}{\bf v}^+\ _{-n}{\bf v}^-$  is a signed circuit of  $_{-A}{\mathcal M}$   and by orthogonality with this circuit  we conclude that  $_{-I}{\bf v}\in X_{-i}^+$. If  $n\in I$  then  orthogonality with the signed circuit  $R({\bf v};i,I)={\bf v}^+\ _{-i}{\bf v}\ _{-Ii}{\bf v}^+ \ _{-I}{\bf v}^-$  also implies that  $_{-I}{\bf v}\in X_{-i}^+$. Therefore  $X_{-i}=(H_i,\emptyset)=Y_{-i}$  is a signed cocircuit of  $_{-A}{\mathcal M}$. 

In a similar way we conclude that $Y_{i}=(H_{-i},\emptyset)$  is also a positive cocircuit of  $_{-A}{\mathcal M}$  and consequently that  $_{-A}{\mathcal M}$  
is an acyclic reorientation of  $\mathcal M$  satisfying condition  1: ${\mathcal F}\subseteq {\mathcal D}'$  of theorem 3.1. Proposition 3.1. then implies that $_{-A}{\mathcal M}$  also satisfies condition 2. ${\mathcal R}\subseteq {\mathcal C}'$  of Theorem 3.1.  
\hskip 5mm

{\parindent=0cm
{\myprop{ Each orientation class of the cube matroid has exactly one orientation that contains  ${\mathcal F}$  as familiy of positive and negative circuits or equivalently that contains  ${\mathcal R}$  as family of signed rank 3 circuits}}\vskip 2mm

{\bf Proof.}} Assume that  $\mathcal M$  and  $_{-A}{\mathcal M}$  are two distinct orientations of  $M(C^n)$  both containing  $\mathcal F$  as family of positive and negative cocircuits.  Then  $_{-A}Y_{n}=(H_{-n}, \emptyset)= \pm Y_n$  and  $_{-A}Y_{-n}=(H_{n}, \emptyset)= \pm Y_{-n}$. 

If $_{-A}Y_{n}= Y_n=(H_{-n}, \emptyset)$  then  $A\cap H_{-n}=\emptyset$  or equivalently  $A\subseteq H_n$. On the other hand,  $_{-A}Y_{-n}=\pm Y_{-n}$  implying that either  $A=\emptyset$  or  $A=H_n$. Since  $\mathcal M \not= _{-A}{\mathcal M}$ it must be  $A=H_n$  but  in this case for all $i\in [n-1]$  the positive circuit  $Y_i=(H_{-i},\emptyset)$  of  $\mathcal M$  verifies  $_{-A}Y_i= (H_{-i}\cap H_{-n}, H_{-i}\cap H_{n})$. Since  both $H_{-i}\cap H_{-n}$  and   $ H_{i}\cap H_{n}$   are nonempty  this contradicts the assumption that  $\mathcal F$  is the subfamily of positive and negative cocircuits of  $_{-A}\mathcal M$.  The case  $_{-A}Y_{n}=- Y_n$  leads to similar contradictions. Therefore  $\mathcal M=_{-A}{\mathcal M}$.  
\vskip 5mm

{\parindent=0cm
\mytheo{ Consider the cube matroid  $M(C^n)$  and the families  ${\mathcal R}$, ${\mathcal F}$, respectively, of the rank 3 signed circuits and  of the positive and negative cocircuits   of the orientation  $Aff(C^n)$  of  $M(C^n)$.

The following three conditions are equivalent:
\begin{enumerate}
\item {\rm (Las Vergnas Cube Conjecture)} $M(C^n)$  has a unique class of orientations.
\item If  ${\mathcal M}$  is an orientation of  $M(C^n)$  containing  ${\mathcal F}$  as signed cocircuits then   ${\mathcal M}=Aff(C^n)$.  
\item  If  ${\mathcal M}$  is an orientation of  $M(C^n)$  containing  ${\mathcal R}$  as signed circuits then   ${\mathcal M}=Aff(C^n)$.  
\end{enumerate}

}
\vskip 2mm

{\bf Proof.}} It is clear from Proposition 3.1. that $2)\Longleftrightarrow 3)$.  The proof that  $2),3)\Longrightarrow 1)$  is a direct consequence of Theorem 3.1. The proof that  $1)\Longrightarrow 2), 3)$  is a direct consequence of Proposition 3.2.   
\vskip 5mm

\section {Final Remarks} 

Theorem 3.2. shows that to prove Las Vergnas Conjecture is equivalent to determine a procedure to 
reconstruct the signature of all the circuits or cocircuits of  $Aff(C^n)$  from the partial subfamilies  $\mathcal F$  and  $\mathcal R$  and the underlying matroid structure. 

We would like to mention that with the description of the n-cube matroid for  $n\leq 7$  in terms of hyperplanes obtained in  \cite {IS1} the (very) interested reader may verify by himself that the signature of cocircuits of  $Aff(C^n)$  can be recovered by orthogonality from  $\mathcal R$  and thus reobtain, in a different way, the result of Bokowski et al \cite {B}.

It is natural to think that if Las Vergnas Conjecture is true then the families  $\mathcal R$  or  $\mathcal F$  might determine not only the orientation  $Aff(C^n)$  of the cube matroid but the oriented matroid  $Aff(C^n)$  itself. The question of whether or not an oriented matroid polytope (the case of  $Aff(C^n)$ ) is determined by its positive cocircuits is known as studying the "rigidity of the matroid polytope" and has been treated in the litterature (see \cite {OM} for a general survey). 
The question of whether or not the family of circuits of fixed rank is enough to determine the oriented matroid has been considered, and studied in a particular case, in \cite {IS2} . 

\vskip 1cm

{

\end{document}